\documentclass[12pt]{article}
\usepackage[utf8]{inputenc} % Umlaute
\usepackage{graphicx} %  LaTeX graphics tool when including figure files
\usepackage{amsmath,amssymb}
\usepackage{color,url}
\newtheorem{proposition}{Proposition}

\newtheorem{remark}{Remark}

\def\re{{\mathbb{R}}}

\def\nz{\hfill\break}

\def\bee{\begin{equation}}
\def\ene{\end{equation}}
\def\beea{\begin{eqnarray}}
\def\enea{\end{eqnarray}}
\def\beeas{\begin{eqnarray*}}
\def\eneas{\end{eqnarray*}}
\def\beas{\begin{eqnarray*}}
\def\eeas{\end{eqnarray*}}

\def\Diag{\hbox{Diag}}
\def\sign{\hbox{sign}}

\def\cN{{\cal{N}}}

\def\cB{{\cal{B}}}
\def\bI{{\mathbf{I}}}
\def\bA{{\mathbf{A}}}
\def\bJ{{\mathbf{J}}}
\def\bK{{\mathbf{K}}}

\def\mr#1{(\ref{#1})}

\def\ignore#1{}

\newenvironment{proof}{\begin{list}{$\!\!${\bf Proof.}%
  \rule{1pt}{0pt}}{\setlength{\leftmargin}{0pt}%
  \setlength{\itemindent}{30pt}%
  \setlength{\listparindent}{15pt}}\item}{\rule{0.3em}{0mm}%
  \hfill\framebox[1.2ex]{\rule{0.3em}{0mm}}\end{list}}

\title { {\bf A Rank-One-Update Method for the
Training of Support Vector Machines}}

\author{
  \small Florian Jarre, Faculty of Mathematics and Natural Sciences,\\
  \small Heinrich Heine Universit{\"a}t D\"usseldorf, Germany
}
\date  {April, 7,  2025}
\begin{document}
\maketitle 

\begin{abstract}
  This paper considers convex quadratic programs
  associated with the training of support vector machines (SVMs).
  Exploiting the special structure of the SVM problem, a new
  type of  active set method with long cycles and stable rank-one-updates
  is proposed and tested (CMU: cycling method with updates).
  The structure of the problem allows for a repeated simple increase
  of the set of inactive constraints while controlling its size. This is 
  followed by minimization steps with cheap updates of a matrix factorization.
  
  A widely used approach for solving SVM problems is the
  alternating direction method SMO,
  a method that is very efficient for generating low accuracy solutions.
  The new active set approach allows for higher accuracy
  results at moderate computational cost. To relate both approaches,
  the effect of the accuracy on the running time and on the
  predictive quality of the SVM is compared based on some numerical examples.
  A surprising result of the numerical examples is that only a
  very small number of cycles (each consisting of less than $2n$ 
  steps) was used for CMU. 
\end{abstract}

\noindent
{\bf Key words:}
Rank-one-update, active set method, support vector machine.

%XXXXXXXXXX github
%findbeta, wmsg damit nur eine warning

\section{Introduction}
An active set descent algorithm is proposed for solving the problem
of training a support vector machine (SVM) 
by exploiting the special structure of the problem and 
solving it with a sequence of {\em cycles}.
Each cycle begins with an  {\em up-cycle} consisting of 
a repeated increase of an {\em inactive set} using very cheap
first-order descent steps. This increase is followed by
a  {\em down-cycle} or 
 {\em sweep} of eliminating  {\em inactive indices} one by one, where
also each Newton step in
the sweep is computationally cheap using rank-one-updates.
Only at the beginning of each sweep a Cholesky factorization
of the  {\em inactive part of the Hessian} is computed.
The overall method is referred to as
cycling method with updates (CMU).

A rank-one-update would also be possible in the up-cycle
while increasing the active set and would save the Cholesky factorization
at the beginning of a sweep. 
However, on the one side, increasing the dimension of the
Cholesky factorization
is more sensitive to numerical rounding errors
than decreasing it. On the other side,
in the up-cycle, the factorization of the Hessian matrix is not needed
for the repeated increases of the active set,
and the overall numerical cost
of a full cycle is comparable to the cost that would have been applicable
when relying on (many)
rank-one-updates in the up-cycle as well. Moreover, the numerical experiments
suggest that for many problems only a very small number of cycles is
necessary, less than 10 in all examples that were tested.
%This numerical observation is interesting also since 
%the primal-dual system of linear programs can be written in a
%similar form.

In many situations 
the method of choice for the training of SVMs
is an alternating direction method referred to as
Sequential Minimal Optimization (SMO), \cite{platt},
that generates approximate solutions
of moderate accuracy in very short time. 

%SMO is reported to generate sparse solutions
%determining $\beta$ by some inactive support vector
%(not at the upper bound $C$)
%for non-trivial problems sweeping is faster

On the other hand, the danger of overfitting the SVM is usually
controlled by the parameters of the kernel and of the so-called
{\em soft margin.}
Solving the SVM problem to high accuracy generally does not
lead to overfitting when the parameters are selected properly.
If the accuracy needed for a given application is not known in advance,
a higher accuracy solution may allow to extract more information 
from a given training set. Some simple numerical examples
that illustrate the benefit of a higher numerical accuracy are
presented in Section \ref{sec.numeric.example}, where the 
active set method of this paper, CMU %(cycling method with updates)
is compared with a variant of SMO
for selected problems, concentrating on the number of arithmetic
operations, final accuracy, and predictive quality.

The setting of SVMs in the present paper follows the outline in
\cite{jarresvm}; for further discussions of SVMs
we refer to \cite{vapnik,schoelbook,HSS,cervantes}
and the references given there.

For SVMs with kernel as considered here, the so-called feature space
(that replaces the data space for the classification task,
see for example, \cite{vapnik,schoelbook,HSS,cervantes,jarresvm}) may
be infinite dimensional, and if not, it is typically high-dimensional.
This implies that for such kernels the low-rank-update formulas that are
the basis of extremely efficient interior-point solvers
for SVMs with low-dimensional data spaces such as
presented in \cite{ferris,masaunders} cannot be applied.
While for SVMs as considered in \cite{ferris,masaunders} the number
of data points may be large, $10^5$ or even $10^6$, the method in
\cite{ferris,masaunders} does not apply to Gaussian or other kernels
as considered in this paper. The method of
the present paper is not suitable for much more than $10^4$ data points.

\subsection{Notation}
$A\succeq 0$ indicates that $A$ is a symmetric positive semidefinite matrix,
and $A\succ 0$ denotes (strict) positive definiteness.
For a vector $x\in\re^n$ the inequality $x\ge 0$ applies componentwise.
%Given a vector $x\in\re^n$, the diagonal matrix with diagonal $x$
%is denoted by $\Diag(x)$ while the diagonal of a square matrix $A$
%is denoted by $\diag(A)$.
Given two vectors $x,y\in\re^{n}$, their Hadamard product
(component-wise product) is denoted by $x\circ y\in\re^n$.
The vector $e:=(1,1,\ldots,1)^T$ always denotes the all-one-vector
with its dimension given by the context. The vector $e_i$ is the
$i$-th canonical unit vector. Let $\bI,\bJ\subset \{1,\ldots,n\}$,
a vector $x\in\re^n$, and a matrix $H\in\re^{n \times n}$ be given.
The vector $x_\bI\in\re^{|\bI|}$
is the vector with components $x_i$ for $i\in\bI$
and $H_{\bI,\bJ}\in\re^{|\bI| \times |\bJ|}$ denotes the submatrix of
$H$ with entries $H_{i,j}$
for $i\in \bI$ and $j\in\bJ$.

\subsection{Basic problem in the training of support vector machines}
\label{sec.problem}
The problem that is to be solved for ``training'' an SVM with kernel and
with soft margins
is a convex quadratic program of the form
\bee\label{qp}
 \min_{x\in\re^n}  \left\{ \frac 12 x^THx
-c^Tx \ \mid \ z^Tx = 0, \ Ce\ge x\ge 0\right\}
\ene
where $H\succ 0$  generally has the form $H=ZKZ$. Here, $K$ is 
a {\em kernel matrix} and $Z$ a diagonal matrix, $C>0$ is a constant
associated with the penalty term of the soft margin,
$z\in\{\pm 1\}^n$ is the vector of {\em labels} for the classification,
and $c=e$ in general. Throughout we assume that the SVM-problem
is well-posed in the sense that $z\not\in \{-e,e\}$. 

The objective function is abbreviated as
$$
q(x)\equiv \frac 12 x^THx-c^Tx,
$$
the box $[0,C]^n$ is denoted by
$$
\cB:= [0,C]^n,
$$
and the null space given by the equality constraint is denoted by
$$
\cN:=\{\ x\mid z^Tx=0\ \}.
$$
Thus, the feasible set of \mr{qp} is $\cB\cap\cN$. 
For later use we note that the orthogonal projection onto $\cN$
is given by $\Pi_\cN:=I-zz^T/n$ because $z^Tz=n$. 
For a feasible point $x$ and some vector $y\in\re^n$,
the projection of $y$ onto the
tangential cone of $\cB$ at $x$ is defined as
$$
p:=\Pi_{\cB,x}(y)
$$
with components
\bee\label{pib}
p_i:=\left\{
\begin{matrix}
  \max \{0,y_i\}  & \phantom{xxxxx} \hbox{if } x_i=0\\
  y_i \phantom{xxxxxxx} &\phantom{xxxxxxxxi}
     \hbox{if } x_i\in (0,C)\\
     \min \{0,y_i\}  & \phantom{xxxxxx}
     \hbox{if } x_i=C.
\end{matrix}\right.
\ene

%Both the evaluation of $\Pi_\cN$ and of $\Pi_{\cB,x}$
%are numerically very cheap. 

\subsection{Optimality conditions}
In Proposition \ref{prop00} below, the 
standard optimality conditions for \mr{qp} are formulated
in a simple form that is used for the
definition of the CMU algorithm.

For a given point $x\in\cB$, the set of
{\em active indices or active inequalities}
at $x$ is always denoted as
\bee\label{ba}
\bA := \{i \mid x_i = 0 \} \cup  \{i \mid x_i = C \}.
\ene
The set of {\em inactive indices} at a given point $x$ is always denoted as
$\bK:=\{1,\ldots,n\}\backslash \bA$.

Denote the gradient of $q$ at $x\in\re^n$
by $g:=g(x):=Hx-c$. For $x\in\cB$ let $\sigma:=\sigma(x)\in\re^n$
be defined as
\bee\label{sigma}
\sigma_i:=\left\{
\begin{matrix}
  1  & \phantom{xxxx} \hbox{if } x_i=0\\
  0  & \phantom{xxxxxxx.} \hbox{if } x_i\in (0,C)\\
 -1  & \phantom{xxxxx} \hbox{if } x_i=C.
\end{matrix}\right.
\ene
Then a direction $s$ at a point $x\in\cB\cap\cN$ is 
a {\em feasible direction} at $x$
(i.e.\  
$x+\lambda s\in \cB\cap\cN$ for small $\lambda >0$) if, and only if,
$z^Ts=0$ and $\sigma\circ s\ge 0$.
The following proposition formulates the optimality conditions for \mr{qp}.

\bigskip

\begin{proposition}\label{prop00}
Let $x\in\cB\cap\cN$, \ $g=Hx-c$,  and
$$
\mu:=\left\{\begin{matrix}
g_\bK^Tz_\bK/|\bK| \phantom{xxxxxxxxxxx} & \hbox{if} \ \ \bK\not=\emptyset \\
\min\{\sigma_ig_i\mid i:\ \sigma_iz_i=1\} &\hbox{else}\phantom{xxxxx}
\end{matrix}\right.
$$
Set $\widetilde g:=g-\mu z$.
Then $x$ is an optimal solution of \mr{qp} if, and only if,
$\sigma\circ\widetilde g\ge 0$ and \, $\widetilde g_\bK=0$.
(Here, $0\in\re^{|\bK|}$ and for $\bK=\emptyset$ the condition
$\widetilde g_\bK=0$ is trivially satisfied.)
\end{proposition}

\bigskip

\begin{proof}
The conditions of Proposition \ref{prop00}
can be rewritten such that they are equivalent to the
standard (necessary and sufficient) KKT conditions.
The proof below is a bit longer, and in case
that the conditions are violated it derives descent directions
that are used in the algorithm. 
  
``$\Rightarrow$'' ($\bK\not=\emptyset$): \ 
Let $x\in\cB\cap\cN$ be an optimal solution of
\mr{qp} with $\bK\not=\emptyset$.
Then the standard KKT conditions imply that
%(since any $s$ with $s_\bA=0$ and with $z_\bK^Ts_\bK=0$
%is a feasible direction at $x$)
there exists $\mu\in\re$ such that
\bee\label{kkt1}
g_\bK=\mu z_\bK,
\ene
i.e.\ $z_kg_k=\mu$ for all $k\in\bK$ (because $z_k^2=1$).
Setting $\mu:=g_\bK^Tz_\bK/|\bK|$ and $\widetilde g:=g-\mu z$,
relation \mr{kkt1} is equivalent to $\widetilde g_\bK=0$.

Now assume that there is an index $i\in\bA$ with $\sigma_i g_i<\sigma_iz_i\mu$.
Then, for some $k\in\bK$ define
\bee\label{descentn}
s:=\sigma_ie_i-\sigma_iz_iz_ke_k.
\ene
It follows that $z^Ts=\sigma_iz_i-\sigma_iz_iz_k^2=0$
and $\sigma\circ s = \sigma_i^2e_i-\sigma_i\sigma_kz_iz_ke_k\ge 0$
since $\sigma_k=0$.
Hence,  $s$ is a feasible direction at $x$.
Moreover, $g^Ts=\sigma_i g_i-\sigma_iz_iz_kg_k
= \sigma_ig_i-\sigma_iz_i\mu
<0$ in contradiction to the optimality of $x$.
Hence, it follows that
\bee\label{kkt2a}
\sigma_ig_i \ge \sigma_iz_i\mu \quad \forall i\in\bA
\qquad \hbox{i.e.} \qquad
\sigma_\bA\circ g_\bA \ge \mu \sigma_\bA\circ z_\bA
\qquad \hbox{or} \qquad\sigma_\bA\circ (g_\bA-\mu z_\bA) \ge 0.
\ene
Since $\sigma_\bK=0$ relation \mr{kkt2a}
is the same as
\bee\label{kkt2}
\sigma\circ (g-\mu z) =\sigma\circ\widetilde g\ge 0.
\ene
``$\Leftarrow$'': 
Conversely, let %\mr{kkt1} and
\mr{kkt2} be satisfied and
a feasible direction $\bar s$ minimizing $g^Ts$ be given,
$$
\bar s\in \hbox{argmin}
\{g^Ts\mid z^Ts=0,\ \ \sigma\circ s \ge 0, \ \ \|s\|_2 \le 1 \}.
$$
Using  \mr{kkt2} and $\widetilde g_\bK=0$,
$$
g^T\bar s=(g-\mu z)^T\bar s
=\widetilde g^T\bar s = (\sigma\circ\widetilde g)^T(\sigma\circ \bar s) \ge 0.
$$
Hence, there is no feasible descent direction starting at
$x$  and thus (by convexity and linearity of the constraints),
$x$ is a minimizer of \mr{qp}. 

\bigskip

``$\Rightarrow$'' ($\bK=\emptyset$): \ 
Similarly, let $x\in\cB\cap\cN$ be an optimal solution
of \mr{qp} with $\bK=\emptyset$.
Then \mr{kkt1}  is void.
Condition \mr{kkt2}, i.e.\ the inequality ``$\sigma\circ \widetilde g\ge 0$'',
is equivalent to %the inequality 
\bee\label{kkt3}
\min\{\sigma_ig_i\mid i: \ \sigma_iz_i=1\}:=\mu\ge
\max\{-\sigma_kg_k\mid k: \ \sigma_kz_k=-1\},
\ene
where the case that $\sigma_iz_i=1$ for all $1\le i\le n$
or $\sigma_kz_k=-1$ for all $1\le k\le n$ cannot occur so that
the min and max are well defined.
(If, for example, $\sigma_iz_i=1$ for all $i$,
then for each $i$, either $\sigma_i=z_i=1$, in which case $x_i=0$,
or $\sigma_i=z_i=-1$,  in which case $x_i=C>0$. And then,
$x^Tz=\sum_{i:z_i=-1} Cz_i<0$ since $z\not\in\{-e,e\}$.
Thus, $x$ is not feasible for \mr{qp}.)

Indeed, when $\bK=\emptyset$ let $i$ be an index with $\sigma_iz_i=1$
and $\sigma_ig_i=\mu$, 
and assume that there is an index $k$
with $\sigma_kz_k=-1$ and $\mu <-\sigma_kg_k$. 
Then 
$$
s:=\sigma_ie_i-\sigma_iz_iz_ke_k = \sigma_ie_i-z_ke_k = \sigma_ie_i+\sigma_ke_k
$$
satisfies $z^Ts=\sigma_iz_i-\sigma_iz_iz_k^2=0$ and $\sigma\circ s\ge 0$.
Hence,  $s$ is a feasible direction at $x$.
Moreover, $g^Ts=\sigma_i g_i-\sigma_iz_iz_kg_k
= \mu-z_kg_k= \mu+\sigma_kg_k
<0$ in contradiction to the optimality of $x$.
Hence, \mr{kkt3} must hold at an optimal solution $x$ with $\bK=\emptyset$.

\noindent
``$\Leftarrow$'': \ 
Conversely, \mr{kkt3} implies $-\mu\le \min\{\sigma_kg_k\mid \sigma_kz_k=-1\}$
i.e.,\ $\sigma_kg_k \ge-\mu$ for $k$ with $\sigma_kz_k=-1$. Using the
definition of $\mu$ in \mr{kkt3} this implies
$\sigma_ig_i\ge \mu\sigma_iz_i$ for all $i\in\{1,\ldots,n\}$
or $\sigma\circ (g-\mu z)\ge 0$.
Given a feasible direction $\bar s$ minimizing $g^Ts$,
$$
\bar s\in \hbox{argmin}
\{g^Ts\mid z^Ts=0,\ \ \sigma\circ s\ge 0, \ \ \|s\|_2 \le 1 \},
$$
and using $\bK=\emptyset$, i.e. $\sigma_i^2=1$ for all $i$,
it follows that
$$
g^T\bar s=(g-\mu z)^T\bar s=((g-\mu z)\circ\sigma)^T(\sigma\circ \bar s)
\ge 0
$$
Hence, there is no feasible descent direction starting at
$x$  and thus, $x$ is a minimizer of \mr{qp}. 
\end{proof}

Summarizing, the necessary and sufficient conditions for optimality
of a feasible point $x$ are also given by \mr{kkt1}, \mr{kkt2}, and \mr{kkt3}.

\section{An active set rank-one-update algorithm}

\subsection{Outline}
From a theoretical point of view, problem \mr{qp} is well understood
and many globally convergent algorithms are available even for more
general convex quadratic programs.
From a practical point of view in turn, the exploitation of the given
structure matters for reducing the overall computational effort,
while generating a solution with sufficiently high accuracy.
In this respect the method of choice in many situations
is the alternating direction method SMO, \cite{platt}.
The present paper is an attempt to improve
over SMO in certain other situations, in particular, when
a numerical solution of high accuracy is needed.

\bigskip

The CMU method proposed in this paper
for solving \mr{qp} always generates feasible iterates. It is divided into
inner iterations and outer iterations:
Each outer iteration consists of a sweep reducing the set of inactive
variables, followed by an up-cycle increasing this set again.
More precisely, each sweep starts at an initial point $x$ in the box
$\cB$ satisfying $z^Tx=0$.
The sweep uses Newton's method to  
successively remove inactive indices one by one using cheap updates
of the Cholesky factor.
When no further inactive indices can be removed, the up-cycle begins,
adding again a moderate number of active indices  to the inactive set.
All steps are such that the objective function decreases.

\bigskip

In more detail, a sweep starting at a point $x$ is as follows:
First, 
a factorization of the inactive part of the Hessian is computed. 
Then, while keeping the active variables fixed at their current value
a Newton step $s$ starting at $x$ is computed for minimizing 
$q$ while maintaining the equality $z^T s=0$
and while changing only the inactive variables.
If the result of the Newton step $x^+:=x+s$ satisfies
$x^+_i \le 0$ or $x^+_i\ge C$, the step length of the
Newton step is reduced so that
$x^+$ lies at the boundary of the box $\cB$,
and the set of inactive variables
is reduced accordingly.
After removing an inactive index, the Cholesky factor is updated.
%(This update is repeated in case that several inactive indices
%reach the boundary of the box $\cB$
%at the same time along the line $\{x+\lambda s\mid \lambda \ge 0\}$.)
Then the  Newton iterations are restarted. This is repeated
until the Newton iterations result in a step that does not lead to a reduction
of the set of inactive indices any more.
A key point of this sweep is that the update of the 
Cholesky factor for the Newton step can be carried out
with order $n^2$ operations for each inactive variable that is removed.
(In general, the Cholesky factor for this application will not be sparse,
but if the Cholesky factor happens to be sparse, there also exist updates that
exploit sparsity, see for example \cite{davis}.)

When no further variable becomes active while performing the Newton step,
the up-cycle begins.
At each step of the up-cycle,
a feasible descent step for $q$ is computed that turns at least one of the
active variables to be inactive.
The step of adding inactive variables
is repeated as long as possible, but only until the active set
has increased by at most 50\%.
This way a moderate number of new inactive variables
are generated before the next sweep is started.
The upper bound of 50\% is intended to keep the size of the Cholesky factor
moderately small.

If the up-cycle fails at the very first attempt to turn one or two of the
active variables to be inactive, the
overall algorithm is terminated.

\bigskip

In what follows the above steps are detailed while observing
the computational cost, numerical accuracy, and global convergence properties.

\subsection{Up-cycle, leaving active constraints}\label{sec.leave}
A key feature of CMU is that in fairly general
situations as detailed in this subsection, iterates can easily be moved
away from parts of the boundary of $\cB$ while decreasing the
objective function $q$.
This allows the repeated increase of the set of inactive
indices before carrying out a numerically expensive recomputation
of a Cholesky factor.

\bigskip

\noindent
The following steps are somewhat ``technical''
but numerically very cheap.

Let a feasible point $x$ for \mr{qp} be given 
and  assume that there is at least one active inequality at $x$,
i.e. $\bA\not=\emptyset$.
In order to reduce the set of active inequalities, i.e.\
to find a feasible point $x^+$ with lower objective function value and
with a smaller set of active inequalities
the following procedure is used.

\bigskip

Define the gradient of $q$ at $x$ as $g:=Hx-c$.
A descent step $\widetilde s$ is then defined as follows.

If the set of inactive indices $\bK$ is not empty,
let $\mu:=g_\bK^Tz_\bK/|\bK|$ and $\widetilde g:=g-\mu z$ 
(see relation \mr{kkt1}). 
Else, let $\mu$ be as defined in
in Proposition \ref{prop00}.
Then set
\bee\label{redefine}
\widetilde g:= g-\mu z \qquad \hbox{and} \qquad
\widetilde s:=\Pi_{\cB,x}(-\widetilde g).
\ene
Then, $\widetilde g^T\widetilde s\le 0$ and for
$\widetilde s\not=0$ the inequality is strict.
Furthermore, $\widetilde s$ is a
feasible direction with respect to all
inequality constraints, i.e. $x+\lambda \widetilde s\in \cB$
for small $\lambda >0$, but in general,
$z^T\widetilde s\not=0$.
The next steps below and in Section \ref{sec.comps} aim at modifying $\widetilde s$ to a
feasible descent direction.
%To this end the role of $\widetilde s$ is detailed next.

Any vector $s\in\cN$
 satisfies  $z^Ts=0$ and
\bee\label{tmp1}
g^Ts=(g-\mu z)^T s=\widetilde g^Ts.
\ene
Therefore $s\in\cN$ is a descent direction for $q$ at $x$ 
whenever $\widetilde g^Ts<0$.

Let $\hat s$ be a feasible  direction for $q$ starting at $x$.
Then, $\hat s$ must satisfy
$\hat s\in\cN$ as well as $\hat s=\Pi_{\cB,x}(\hat s)$.
And using \mr{tmp1} it follows that
\bee\label{opt1}
g^T\hat s=\widetilde g^T\hat s 
\ge  -(\Pi_{\cB,x}(-\widetilde g))^T \hat s =
-\widetilde s^T\hat s,
\ene
where the inequality follows from relation \mr{technical} below:

\begin{itemize}
\item
Let $x$ be feasible for \mr{qp} and 
let $\hat s$ be given with $\hat s=\Pi_{\cB,x}(\hat s)$. Then, for arbitrary $g$
it follows that
\bee\label{technical}
g^T\hat s \ge -(\Pi_{\cB,x}(-g))^T\hat s.
\ene
Indeed, let $I:=\{i\mid x_i=0\}$,
$J:=\{j\mid x_j\in (0,C)\}$, and $K:=\{k\mid x_k=C\}$.
Then,  $\hat s_I\ge 0$
and $\hat s_K\le 0$. Let $p:= \Pi_{\cB,x}(-g)$, then
$p_I=\max(0,-g_I) \ge -g_I$ and $p_K=\min(0,-g_K)\le -g_K$
where min and max are taken componentwise, and
$$
-g^T\hat s= \underbrace{-g_I^T}_{\le p_I^T}\underbrace{\hat s_I}_{\ge 0}+
\underbrace{-g_J^T}_{= p_J^T}\hat s_J+
\underbrace{-g_K^T}_{\ge p_K^T}\underbrace{\hat s_K}_{\le 0}
\le p_I^T\hat s_I+p_J^T\hat s_J+p_K^T\hat s_K = p^T\hat s.
$$
Therefore,
$$
g^T\hat s \ge -p^T\hat s = -(\Pi_{\cB,x}(-g))^T\hat s. 
$$
{\rule{0.3em}{0mm}  \hfill\framebox[1.2ex]{\rule{0.3em}{0mm}}}
\end{itemize}

\bigskip

\noindent
Hence, any feasible descent direction  $\hat s$ starting at 
$x\in\cB\cap\cN$  satisfies $0>g^T\hat s\ge -\widetilde s^T\hat s$.
In particular, the following proposition is true:

\begin{proposition}\label{prop1} 
  When $x$ is feasible for \mr{qp}, $g:=Hx-c$, and $\widetilde s$
  given by \mr{redefine}  satisfies $\widetilde s=  0$ 
  then there does not exist a feasible descent direction,
  i.e. $x$ solves the convex problem \mr{qp}. 
\end{proposition}

\subsubsection{\bf Computing a search direction:} \label{sec.comps}
In the following, for a feasible point $x$ let $\widetilde s\not = 0$
be as in \mr{redefine} and let
\bee\label{bI}
\bI:=\{i\in\bA \mid z_i\widetilde s_i > 0\}\quad \hbox{and} \quad
\bJ:=\{i\in\bA\mid z_i \widetilde s_i < 0\}
\ene
be the sets of active indices that increase/decrease the constraint
term $z^T \widetilde s$.
%when moving along $\lambda\widetilde s$.
For defining a feasible descent direction that makes certain
active variables inactive, the following cases are considered:

\bigskip

\noindent
{\bf Case 1:}
    
If $\bI$ and $\bJ$ are both nonempty, define a feasible direction
$s$ as follows:
\bee\label{iandj}
v_1 := z_\bI^T\widetilde s_\bI > 0, \quad v_2 := z_\bJ^T\widetilde s_\bJ < 0,\quad
s_\bI := -v_2 \widetilde s_\bI , \quad s_\bJ := v_1 \widetilde s_\bJ ,
\ene
and $s_i := 0$ for all $i$ not in $\bI\cup\bJ$.
Then,
$$
z^Ts=z_\bI^T s_\bI+ z_\bJ^T s_\bJ
= -v_2z_\bI^T \widetilde s_\bI+ v_1 z_\bJ^T\widetilde s_\bJ= -v_2v_1 + v_1v_2
=0
$$
and $\widetilde g^T s <0$ since all nonzero components of $s$
have opposite sign of the associated components in $\widetilde g$. Moreover,
by the sign structure of $s$, it is a feasible descent direction
that has at least two nonzero active components.
So, at least two active components become inactive along the line
$x+\lambda s$ for $\lambda >0$, and large values of $\lambda$
are possible without violating any inequality constraint since
$s$ only moves active indices ``away from the boundary''
(see also the last paragraph in Section \ref{sec.technical}).

\bigskip

\noindent
{\bf Case 2:}
   
When exactly one of the sets  $\bI$ or $\bJ$ is empty
let $i\in\bA$ be an active index with maximal value of $|\widetilde s_i|$.
By definition of $\bI$ and $\bJ$, it follows
that $|\widetilde s_i|>0$ and $\sign(\widetilde s_i) = \sigma_i =
-\sign(\widetilde g_i)$.

Let
\bee\label{defbk}
\overline\bK:=\{j\mid \sigma_i\sigma_j z_iz_j \le 0\}.
\ene
Since $\sigma_j=0$ for $j\in\bK$, it follows that  $\bK\subset\overline\bK$.
In addition, $\overline\bK$ contains active indices $j$ such that
 the vector
\bee\label{iorj}
s:=\sign(\widetilde s_i)(e_i-z_iz_je_j)
\ene
is a feasible direction with respect to all constraints of \mr{qp}.
Indeed, $z^Ts=\sign(\widetilde s_i)(z_i-z_iz_j^2)=0$ and
$\sigma\circ s=\sigma_i^2e_i-\sigma_i\sigma_j z_iz_je_j\ge 0$.
\ignore{
as shown next, $x+\lambda s\in\cB$ for small $\lambda>0$:
The $i$-th component of $x+\lambda s$ is given by
$x_i+\lambda \,\sign(\widetilde s_i)\in (0,C)$ for $\lambda \in (0,C)$
by definition of $\widetilde s$.
And either $j$ is inactive, remaining inactive for small $\lambda>0$,
or for small $\lambda >0$, the $j$-th component of $x+\lambda s$ is given by
$$
x_j-\lambda \,\sign(\widetilde s_i)z_iz_j=
x_j-\lambda \,\sigma_iz_iz_j=
\left\{\begin{matrix}
C+\lambda\sigma_i\sigma_j z_iz_j < C& \hbox{if} \ \sigma_j = -1 \\
0-\lambda\sigma_i\sigma_j z_iz_j > \, 0& \hbox{if} \ \sigma_j = 1,
\phantom{c} \end{matrix}\right.
$$
where the inequality follows from the definition of $\overline\bK$.}

The linearization of the objective function along the above
search direction is
$$
g^Ts = \widetilde g^Ts = \widetilde g^T\sign(\widetilde s_i)(e_i-z_iz_je_j)
=\widetilde g_i\sign(\widetilde s_i)-\sign(\widetilde s_i)z_iz_j\widetilde g_j
= -|\widetilde g_i| -\sign(\widetilde s_i)z_iz_j\widetilde g_j,
$$
since $\sign(\widetilde g_i)=-\sign(\widetilde s_i)$.
Therefore an inactive index $j$ is selected with maximum value of
$\sign(\widetilde s_i)z_iz_j\widetilde g_j$ in order to define $s$.
The maximum value may be negative, and if it is less or equal to 
$-|\widetilde g_i|$, the resulting $s$ is not a descent direction.
In this case the increase of the inactive set fails and the up-cycle
is terminated.

\bigskip 

\noindent
{\bf Case 3:}
If $\bI=\bJ=\emptyset$ then the increase of the inactive set fails as well
and the up-cycle is terminated.

\bigskip

\begin{remark}
Assume that the up-cycle fails at the first step
immediately after the sweep with Newton iterations is completed.
Let $x$ be the iterate generated by the last Newton iteration 
and denote the sets of active and inactive indices at $x$ by $\bA$ and $\bK$,
respectively.
By Newton's method, either $\bK=\emptyset$ or $\widetilde g_\bK = 0$.
Therefore, if failure happens in Case 3, then 
$\widetilde s = 0$, and  by Proposition \ref{prop1},
$x$ is an optimal solution of \mr{qp}.

If $\bK$ is not empty, then,
since $\widetilde g_\bK=0$, the failure
in the up-cycle cannot happen in Case~2. Hence in this case as well,
$x$ is an optimal solution of \mr{qp}.
When  $\bK=\emptyset$, i.e., when $\bA=\{1,\ldots,n\}$
failure could happen in Case~2.
In this case, by definition of $i$, the value $\sigma_ig_i$ either
coincides with the ``min'' in \mr{kkt3} or with the ``max'', and by the
criterion that leads to the failure, the inequality in \mr{kkt3}
is violated. Again, it follows that $x$ is a minimizer of \mr{qp}.
\end{remark}

\subsubsection{Line search}\label{sec.linesearch}
If a feasible descent direction $s$ was found in Case 1 or Case 2,
a  line search is carried out along the direction $s$
minimizing $q(x+\lambda s)$ for $\lambda \ge 0$, i.e.\
first define the exact line search step length
$$
\hat \lambda := -\nabla q(x)^Ts/(s^T Hs)>0
$$
and the maximum feasible step length along $\lambda s$,
$$
\lambda_{max} := \min\left\{\  
\min_{i:s_i>0}\{(C-x_i)/s_i\}, \ \ \min_{i:s_i<0}\{-x_i/s_i\} \right\}.
$$
Then let $\lambda := \min\left\{\ \hat \lambda, \ \ %\frac 12
\lambda_{max} \right\}$ and define the point 
\bee\label{lines}
x^+:= x+\lambda s.
\ene
%(The factor $\frac 12$ prevents that another inquality becomes active.)
If $\lambda = \lambda_{max}$ a new active constraint is added
while at least one other active constraint becomes inactive.
In principle, this might possibly lead to cycling leaving and adding the
same constraints along shorter and shorter steps.
Therefore the up-cycle is limited to at most $n$ steps
before starting another sweep.

Summarizing, the up-cycle is as follows:

\subsubsection{Summary, up-cycle}
Let $\bA\not=\emptyset$ be as in \mr{ba} and define
$\widetilde s$ as in \mr{redefine}.\nz
Set max$\bA := \min\{n, \max\{ 100,\lceil 3|A|/2\rceil\}\}$, $k:=0$.\nz
Set $\bI:=\{i\in\bA\mid z_i\widetilde s_i > 0\}$ and $\bJ:=\{i\in\bA\mid z_i
\widetilde s_i < 0\}$.\nz
Repeat Step 1.\ -- Step 3.\ until a Stop command is encountered:
\begin{enumerate}
\item
  {\bf Case 1:} If $\bI$ and $\bJ$ are both nonempty,
  define $s$ as in \mr{iandj}.

  {\bf Case 2:}
  If exactly one of the sets $\bI$ or $\bJ$ is nonempty,
  let $i$ maximize $|\widetilde s_i|$ for $i$ in $\bA$.
  Then $j$ is selected as to maximize
  $\sign(\widetilde s_i)z_iz_j\widetilde g_j$
  for all indices $j$ in $\overline\bK$ defined in \mr{defbk}.
  Then define $s$   as in \mr{iorj}.

  {\bf Case 3:}
  If both the sets $\bI$ or $\bJ$ are empty, set $s=0$.
\item
  If $s$ generated in Step 1. is not a (strict) descent step,
  Stop. The further reduction of the set of  active indices fails.

  Else, %do a line search along $x+\lambda s$ and 
  define the next iterate $x:=x^+$ with $x^+$ given in \mr{lines}.
  Update $\bA$, $\bI$, $\bJ$, and 
  set $k:=k+1$. If $k\ge n$ or if $|A| \ge$ max$\bA$, Stop.
\end{enumerate}

\begin{remark}
If the step length in the up-cycle leads to the
boundary of $\cB$, then the cardinality of $\bA$ might not increase
for this step. To avoid discussions of possible cycling when this case
occurs repeatedly,
the safeguard query ``If $k\ge n$'' is added in Step 3.\ above.
(And the restriction, not to force max$\bA$ below 100 is subject to change;
it is merely intended to reduce the number of cycles.) 
\end{remark}

\subsection{Starting point}\label{sec.start}

To find a suitable starting point for the overall algorithm,
the following procedure is used.
Set $s:=\Pi_\cN c$.
As $c=e$ and $z\in \{\pm 1\}^n \backslash \{-e,e\}$
it follows that  $|z^Tc| < z^Tz$ and $s > 0$.
Then do a line search as in Section \ref{sec.linesearch} minimizing
$q$ along $q(0+\lambda  s)$. When $\lambda < \lambda_{max}$
the result of the line search
defines a point where all variables are inactive.
In this case, when $n$ is large, 
the computational effort for computing a Cholesky
factor at this point will be large. To limit the cost of
the Cholesky factor, another point is chosen as starting point for CMU:
%The above point $x=\lambda  s$
%will not be used as a starting point.
From the above point $\lambda s$ determine a set of indices
of moderate cardinality that are
``promising'' (with the largest values $x_i-(\Pi_\cN\nabla q(x))_i$)
and define a search step $s$ using only these promising  coordinates
projected onto $\{s\mid z^Ts=0\}$. Then
repeat the line search as detailed above.

\subsection{Down-cycle, Newton iterations}\label{sec.newt}
Let a feasible point $x$ with set of active indices
$\bA$ be given.
Let $\bK$ be the set of the remaining (i.e. the inactive) indices.
Keeping $x_\bA$ fixed in $x =(x_\bA^T,x_\bK^T
+ \Delta x_\bK)^T$,
the minimization of $q$ with respect to $\Delta x_\bK$
subject to $z_\bK^T\Delta x_\bK =0$ is considered next.
As described in Section \ref{sec.problem}, the objective function is given by
$$
q(x) = \frac 12 x^TH x - c^Tx
= \frac 12(x_\bA^T,x_\bK^T
+ \Delta x_\bK^T) H (x_\bA^T,x_\bK^T
+ \Delta x_\bK^T)^T  - c^T(x_\bA^T,x_\bK^T
+ \Delta x_\bK^T)^T
$$
$$
=\frac 12 \Delta x_\bK^TH_{\bK,\bK}\Delta x_\bK -
(c_\bK-H_{\bK,\bA}x_\bA- H_{\bK,\bK}x_\bK)^T
\Delta x_\bK + const
$$
$$
=\frac 12 \Delta x_\bK^TH_{\bK,\bK} \Delta x_\bK -
(c-Hx)_\bK^T
\Delta x_\bK + const,
$$
where $const$ is a term that does not depend on $\Delta x_\bK$.
The Newton step for minimizing $q$ within the set $\cN$ is given by
the linear system
$$
H_{\bK,\bK} \Delta x_\bK + \eta z_\bK = (c-Hx)_\bK, \qquad z_\bK^T\Delta x_\bK=0. 
$$
where $\eta\in\re$ is a Lagrange multiplier.
Let
$$
[u,v]:=H_{\bK,\bK}^{-1}[(c-Hx)_\bK,z_\bK]
$$
be evaluated by using a Cholesky factorization of $H_{\bK,\bK}$.
By positive definiteness of $H_{\bK,\bK}$, it follows $z_\bK^Tv>0$ and 
$\eta:=(z_\bK^Tu)/(z_\bK^Tv)$ is well defined. Set
\bee\label{newt}
\Delta x_\bK:=u-\eta v.
\ene
Then, $z_\bK^T\Delta x_\bK=z_\bK^T u- \tfrac{z_\bK^Tu}{z_\bK^Tv} z_\bK^Tv=0$
and
$$
H_{\bK,\bK} \Delta x_\bK + \eta z_\bK =H_{\bK,\bK}(u-\eta v)+ \eta z_\bK
=(c-Hx)_\bK,
$$
i.e.\ \mr{newt} defines the Newton step.

The next iterate is defined using the maximum step length $\le 1$
that maintains the bound constraints. 
If a new constraint becomes active at the next iterate,
the sets $\bA$, $\bK$ and the Cholesky factor are updated
with $O(n^2)$ operations as outlined below.

Else, if no new constraint becomes active at the next iterate,
the sweep terminated. In this case
$x_\bK$ is a global minimizer
of $q$ with respect to indices in $\bK$ and 
subject to the equality constraint $z_\bK^Tx_\bK = 0$, i.e.
$\exists \mu \in\re: \ g_\bK = \mu z_\bK$.

\bigskip

\subsection{Cholesky updates}
There are five different rank-one-updates for a Cholesky factor
presented in \cite{ggms} [p.514-523]. A rank-one-update also
is available in {\tt cholupdate}  in Matlab
and this is used in the numerical examples below.
Octave also offers a command {\tt choldelete}
performing an elimination of a row and column
of the matrix to be factored. 
Following the outline in \cite{pss}
it is briefly described next how the deletion of a row and column
amounts to a rank-one-update as available in {\tt cholupdate}.

If a Cholesky factorization is given,
$$
LL^T=
\left[ \begin{matrix}
L_{1,1} & 0 & 0 \\ l_{1,2}^T & l_{2,2} & 0 \\ L_{3,1} & l_{3,2} & L_{3,3}
\end{matrix}\right]
\left[ \begin{matrix}
L_{1,1}^T & l_{1,2} & L_{3,1}^T \\ 0 & l_{2,2} & l_{3,2}^T \\ 0 & 0 & L_{3,3}^T
  \end{matrix}\right]
=
\left[ \begin{matrix}
    A_{1,1} & a_{1,2} & A_{3,1,}^T \\ a_{1,2}^T & a_{2,2} & a_{3,2}^T
    \\ A_{3,1} & a_{3,2} & A_{3,3}
\end{matrix}\right]
$$
with lower case letters indicating column vectors and upper case letters
indicating (sub-) matrices of appropriate dimensions,
and if the row $( a_{1,2}^T , a_{2,2} , a_{3,2}^T)$ and the associated
column are to be deleted such that 
$$
\widehat L\widehat L^T=
\left[ \begin{matrix}
    \widehat L_{1,1} &  0 \\
    \widehat L_{3,1} &  \widehat L_{3,3}
\end{matrix}\right]
\left[ \begin{matrix}
    \widehat L_{1,1}^T & \widehat L_{3,1}^T \\
    0 &  \widehat L_{3,3}^T
  \end{matrix}\right]
=
\left[ \begin{matrix}
    A_{1,1}  & A_{3,1,}^T 
    \\ A_{3,1}  & A_{3,3}
\end{matrix}\right]
$$
then $\widehat L_{1,1}\widehat L_{1,1}^T=A_{1,1}=L_{1,1}L_{1,1}^T$  implies
$\widehat L_{1,1}=L_{1,1}$. From
$\widehat L_{3,1}\widehat L_{1,1}^T=A_{3,1}=L_{3,1}L_{1,1}^T$ then follows
$\widehat L_{3,1}=L_{3,1}$. Finally,
$$
L_{3,1}L_{3,1}^T + l_{3,2}l_{3,2}^T + L_{3,3}L_{3,3}^T = A_{3,3} =
\widehat L_{3,1}\widehat L_{3,1}^T+\widehat L_{3,3}\widehat L_{3,3}^T=
 L_{3,1} L_{3,1}^T+\widehat L_{3,3}\widehat L_{3,3}^T
$$
implies $\widehat L_{3,3}\widehat L_{3,3}^T = L_{3,3}L_{3,3}^T+ l_{3,2}l_{3,2}^T$,
a rank-one-update of $L_{3,3}$.
This update is stable as the term $ l_{3,2}l_{3,2}^T$
is added to the existing Cholesky factor.
(Augmenting a given factorization by a row and a column is
numerically less stable, in general.)

\subsection{Global convergence}
Global convergence of the algorithm follows immediately
from the derivations of the previous sections using a standard argument of
active set methods: At the end of each sweep, the minimizer of $q$
on the active set is generated. Since each step of the algorithm
is a (``strict'') descent step, this active set will never be revisited again.
As there are only finitely many different active sets, the algorithm
will terminate after a finite number of outer iterations,
each of which taking at most
$n$ steps in the up-cycle and at most $n$ steps during the sweep
following thereafter.

\subsection{Technical details}\label{sec.technical}

%When a sweep terminates, i.e.  
%when no further variable becomes active while performing the Newton step,
%one step of iterative refinement for this Newton step is applied
%to increase the numerical accuracy, and then the up-cycle begins.
%By refactoring at the end of each sweep, the accumulation of rounding errors
%during the preceeding Cholesky updates steps will be reduced. 

When performing a line search along a Newton direction,
it might happen that two inactive indices turn active at the same time.
For simplicity only one index is added to the active set.
In the next Newton step, the step length might be zero and the set of
active indices will then be increased again using another update
of the Hessian matrix.  

For Gaussian kernels with small exponents, the matrix $H$ may be
very ill-conditioned and the norm of the optimal solution of \mr{qp} may
be huge. In such situations, interior-point methods generally are
unstable. Large active upper bounds such as $C=10^{10}$ may lead to a rather
high numerical error in the optimal solution. (In fact, this was the reason
why the work on the present paper using active sets was started.)
For the numerical experiments in Section
\ref{sec.numeric.example}, a small multiple of the identity was added
(as a regularization term) to 
$H$ before forming the Cholesky factor. The result of the regularized Newton
step then was corrected using one or two steps of iterative refinement
(where the residuals for the iterative refinement
corrections were computed with the original matrix $H$.)
When the norm of the optimal solution often is huge
the computation of $g_\bK = (Hx-c)_\bK$
frequently is subject to large cancellation errors.  
To take this into account a relative KKT condition is
listed in the numerical section, using $\|\widetilde g_\bK\|_2/\|x\|_\infty$
in place of $\|\widetilde g_\bK\|_2$. In such situations also the
evaluation of $q(x)$ may be subject to high cancellation errors.

As an obvious technical detail, 
in the numerical implementation the definition \mr{pib} of $p$
is changed using an $\epsilon$-tolerance: more precisely,
the cases $x_i\le \epsilon$,  $x_i\in (\epsilon,(1-\epsilon)C)$
and $x_i \ge (1-\epsilon)C$ are distinguished
where  $1 \gg \epsilon >0$ is a user-defined tolerance for the active set.
Likewise, definition \mr{ba} is modified
$\bA := \{i \mid x_i \le \epsilon \} \cup  \{i \mid x_i \ge (1-\epsilon)C \}$,
and same for the definition \mr{sigma}.

We close the discussion of technical details with the remark that for
small positive values of $C$, a step in the up-cycle might lead to a step length
that makes another inactive variable active, and thus leads to
repeated short steps in the up-cycle, oscillating between lower and upper bound.
Since typical values of $C$ for SVM are rather large, 
we did not concentrate on this case but just stopped the up-cycle
after at most $n$ steps.
Alternatively, one could  shorten the step length in the up-cycle to ensure
that no
inactive variable becomes active.

\section{A greedy SMO Algorithm}
The SMO algorithm \cite{platt} consists of repeating 
the choice of $s = \pm e_i \pm e_j$ for properly selected indices $i\not=j$
as in Step 2.\ of the up-cycle, followed by a line search
$x\leftarrow x+\lambda s$.
The rules used in the up-cycle of CMU, however, differ from the
rules for SMO since CMU selects
active indices $i,j$ whenever possible.

It is well known, see \cite{powell} for example, that
alternating direction methods such as SMO may fail to converge if the
selection of the directions is not carried out carefully,
an aspect that is addressed later in this section.
The original ``SMO-paper'' \cite{platt} refers to a general class of
algorithms to establish global convergence.
In several subsequent papers, see for example \cite{keerthi}
and the references therein,
numerous variants of the selection of $i$ and $j$ in the SMO algorithm have
been proposed
and their convergence has been established.

For the numerical experiments in this paper
a greedy selection of the indices $i,j$ in the
SMO algorithm as well as a randomized
selection are compared. The computational cost of the greedy selection is
of order $n$ operations, i.e.\ of the same order as the cost for
one SMO step with any other pivoting rule, and it is substantially cheaper
than one step of the sweeping cycle of CMU.

Before addressing the shortfalls of the greedy selection  
it is outlined next:

At each step of the greedy SMO,
the projected gradient $\widetilde g$ (see \mr{redefine})
is used and updated.
Let a feasible iterate $x$ be given and set $g:=Hx-c$ and
$\widetilde g:=\Pi_\cN g=g-\mu z$ with $\mu=g^Tz/n$.
(This choice of $\mu$ differs from Proposition \ref{prop00}.)

A sparse (numerically cheap) search
direction $s$ starting at a feasible iterate $x$ of \mr{qp}
is defined using two indices $i$ and $j$ and setting
$\widetilde s:=\Pi_{\cB,x}(-\widetilde g)$  as in \mr{redefine}.
Then, $s :=\sign(\widetilde s_i)(e_i-z_iz_je_j)$. 
As in \mr{iorj} it follows that $z^Ts=0$
and
$$
q(x+\lambda s) = q(x) + \lambda \sign(\widetilde s_i)(g_i-z_iz_jg_j) +
\tfrac{\lambda^2}2 (H_{i,i}+H_{j,j}-2z_iz_jH_{i,j}).
$$
Using \mr{tmp1} the relation $z^Ts=0$ implies
  $(g_i-z_iz_jg_j)=g^Ts=\widetilde g^Ts
=(\widetilde g_i-z_iz_j\widetilde g_j)$  and
\bee\label{qij}
q(x+\lambda s) = q(x) +
\lambda   \sign(\widetilde s_i)   (\widetilde g_i-z_iz_j\widetilde g_j) +
\tfrac{\lambda^2}2 (H_{i,i}+H_{j,j}-2z_iz_jH_{i,j})
\ene
where $\|\widetilde g\|_2\le\|g\|_2$.
%and $\sigma_i=\sign(\widetilde s_i)=- \sign(\widetilde g_i)$.
Based on \mr{qij}, the indices $i$ and $j$ are chosen successively:
\begin{enumerate}
  
\item
First, a greedy approach  selects $i$
as the entry maximizing $|\widetilde s_i|$.
Then set 
$$
\overline\lambda_i :=
\left\{\begin{matrix}
\phantom{C  } -x_i 
& \hbox{if} \  \ \widetilde s_i \le 0   \phantom{.}\\
 C-x_i  & \hbox{if} \ \ \widetilde s_i > 0 .
\end{matrix}\right.
$$

\item
Given $i$ and the associated value $\overline\lambda_i$,
again a greedy heuristics is used to define $j$ such that
$q$ in \mr{qij} is minimized along $s=\sign(\widetilde s_i)(e_i-z_iz_je_j)$
subject to box constraints.
More precisely, for $j\not = i$ with
$\widetilde g_i(\widetilde g_i- z_i z_j\widetilde g_j)> 0$
let\footnote{%
The restriction $\widetilde g_i(\widetilde g_i- z_i z_j\widetilde g_j)> 0$
implies that $\sign(\hat\lambda_j)=\sign(\overline \lambda_i)
\ \ (=\sign(\widetilde s_i)=-\sign(\widetilde g_i))$.
}
$$
\hat\lambda_j:=
-(\widetilde g_i- z_i z_j\widetilde g_j)/(H_{i,i}+H_{j,j}-2z_iz_jH_{i,j}).
$$
Then $\overline\lambda_j$ is determined as
$$
\overline\lambda_j := \left\{\begin{matrix}
\max\{  \hat\lambda_j,\phantom{C  }  -x_j \}
& \hbox{if} \ \ \hat\lambda_j \le 0 \phantom{.}\\
\min\,\{  \hat\lambda_j, C-x_j \} & \hbox{if} \ \ \hat\lambda_j > 0 .
\end{matrix}\right.
$$
Finally, set $\overline\lambda_j:=\sign(\overline\lambda_j)
\min\{ |\overline\lambda_i|,|\overline\lambda_j| \} $.
Then $j$ is selected as to minimize
$$
\overline\lambda_j (\widetilde g_i-z_iz_j\widetilde g_j) +
\tfrac{\overline\lambda_j^2}2 (H_{i,i}+H_{j,j}-2z_iz_jH_{i,j}).
$$
\end{enumerate}

\begin{remark}
For Gaussian kernels satisfying $H_{i,i}\equiv 1$ for all $i$,
the selection of $j$ can be carried out with about $20n$ floating point
operations. For large values of $n$ this effort may pay off by
the reduction of the objective value.
Choosing $i$ and $j$ simultaneously might result in an
even larger reduction
of $q$ but generally, this would require order $n^2$ operations considering
all entries $H_{i,j}$ for $1\le i\le j\le n$.
\end{remark}

\begin{remark}
The heuristics of choosing $i$ can result in a choice of $i$
with $x_i$ close to the boundary of $\cB$ and such that 
only a very short step length is possible, no matter how $j$ is chosen.
In this case, however, $i$ will become active in the next step and $\tilde s_i$
will then be zero.
For several other seemingly profitable choices of $i$
one can construct examples that may lead to convergence to non-optimal points.
\end{remark}

When $i,j$ have been selected, the line search along $s=
\sign(\widetilde s_i)(e_i-z_iz_je_j)$
leads to the step length $\overline\lambda_j$ (possibly negative)
that was computed during
the selection process of $i,j$.

Then $g$ is updated as $g\leftarrow g+ \overline\lambda_j(He_i-Hz_iz_je_j)$
with $3n$ floating point operations.
Thereafter $\widetilde g$ is updated as $\widetilde g\leftarrow g-
\frac{g^Tz}n z$ with another $3n$ operations.

The greedy SMO algorithm stops when
$\mr{kkt1}, \mr{kkt2}, \mr{kkt3}$ are satisfied up to 
some tolerance $\epsilon$ or when a given maximum number of
iterations has been reached.

Choosing $i$ and $j$ uniformly randomly from $\{1,\ldots,n\}$
is about 10 times cheaper than the above greedy heuristic,
and it is guaranteed to converge without the danger of running into a
cycle. However, in the numerical examples, the random choice
is more than 10 times slower than the greedy heuristic; the latter one therefore
is used for a conceptual comparison.

\section{Numerical experiments}\label{sec.numeric.example}
Some numerical experiments were carried out to compare the overall solution
times and the final accuracy of the solutions generated by CMU.
The experiments were carried out using Gaussian kernels.
As argued in \cite{jarresvm}, the Gaussian kernels are optimal
with respect to a self-concordance parameter similar to the one
introduced for barrier functions in \cite{nene}.

Two sets of examples were used, half-moon shapes and checker board patterns.
The solution times are always listed in seconds.

\subsection{Half-moon shapes} \label{sec.half}
The first set of examples uses a higher-dimensional half-moon shape,
i.e. a connected non-convex set in $d\ge 2$ dimensions.
% generate n sample points x^i \in R^d for a support vector machine:
The input for this example is: $d\ge 2$ (dimension of the data space),
$\delta \in (0,2)$,
%$\rho\ge \delta$,
and $n$ (number of training points).

In the numerical experiments, unless stated otherwise, the parameters
$$\delta=\frac 14\qquad \hbox{and} \qquad n=500$$
are used. 

First, define a set 
$$
S_1 = \{ x \in \re^d \ \ |  \ \ \|x\|_2 \le 1, \ \|x+\delta e_1\|_2 \ge 1 \}.
$$
For $\delta=\frac 14$ and for any $i\in
\{2,\ldots,d\}$ the projection of $S_1$ onto the $(x_1,x_i)$-plane 
has a half-moon shape. In general dimensions $d\ge 2$ the set $S_1$
is the difference
of two Euclidean balls that have non-empty intersection.

The SVM then is to decide whether a given "new"
data point lies in $S_1$ or not.

By definition, $S_1 \subset [-1,1]^d$ but for large $d$ the volumes are
vol($S_1$) $\ll$  vol($[-1,1]^d$) = $2^d$ even when $\delta$ is large.
In this case, drawing the training points uniformly from $[-1,1]^d$
would result in extremely unbalanced labels.
Likewise for the test points.

Therefore define $S^- = S_1 - \delta e_1$.  
For $x \in S^-$, it then follows $x + \delta e_1  \in S_1$, 
i.e.\ $\| x+\delta e_1 \|_2 \le 1$, so that the interiors of $S_1$ and $S^-$
are  disjoint.
Similarly let $S^+ = S_1 + \delta e_1$; then also the interior of
$S^+$ is disjoint from $S_1$. Then set 
$S_2 = S^- \cup S^+$.

The sample points $x=x^i$ are then generated with some random
distribution\footnote{%
More precisely, a point $x$ in $S_1$ is generated as follows:
First $x_1$ is generated uniformly from $[-\delta/2,1]$.
Then, a random normal vector $y \in \re^{d-1}$ is drawn and
$\bar y := y / \|y\|_2$ is defined.
Then, $x := [x_1; \lambda y]$ where $\lambda$ is drawn uniformly from
$[\delta_1, \delta_2]$ with 
$\delta_1 = \sqrt{\max\{0, 1-(x_1+\delta)^2\}}$,
$\delta_2 = \sqrt{1-x_1^2}$.
Half of the training points are then shifted to both parts of
the set $S_2$. The training points generated
this way are concentrated at both ends of the
``half moon'' and fewer points in the middle, an effect that is
even more pronounced for $d>2$.
}
within the union of $S_1$ and $S_2$ and with labels 
$z_i =  1$ if $x \in S_1$ and $z_i = -1$ if $x \in S_2$.
(The definition of $S_2$ implies that the SVM is to
generate a two-sided approximation of $S_1$.)

\ignore{
\vfill
\eject
Figure 1 illustrates the location of the training points generated
when $d=2$, with training points concentrated at both ends of the
``half moon'' and fewer points in the middle, an effect that is
even more pronounced for $d>2$.
\vbox{
  \vskip-4cm
%\centerline{
\hskip-5cm  \includegraphics[width = 12cm]{./half_dist}
%}
\vskip-5cm
%\centerline{
\hskip-5cm   {\bf Figure 1}, Distribution of 2000 training points for $d=2$.
%}
\bigskip
\bigskip
}
}

For the numerical examples below, $10^5$ test points drawn from
the same distribution were always used
to assess the classification error.

\subsubsection{Parameter selection}\label{sec.param}
For simplicity, all training points were classified correctly and
the soft margin constant $C$ therefore was set to $C=\infty$.
First a comparison of the constants $\gamma$ in the Gaussian kernel
was carried out for $d=2$ dimensions of the data space.
(The kernel is based on the function $e^{\gamma\|x-y\|^2_2}$.)
%and follows the setting given in \cite{jarresvm}.)
The results are shown in Table~1.

\begin{center}
\begin{tabular}{ c c c c c c c c c }
  $\gamma$ & cycles & iterations & time & KKT violation &
  $q(x^{final})$ & $\|x^{final}\|_\infty$ & rel class. errors\\
  0.03 & 3 & 433  & 0.69 & 1.8e-11 & -1.0e+12 & 5.8e+10 & 0.0106, \ \ 0.0267 \\
  0.3  & 7 & 1328 & 1.21 & 4.3e-16 & -3.4e+09 & 2.9e+09 & 0.0327, \ \ 0.0223 \\
  3    & 6 & 832  & 0.85 & 5.1e-16 & -4.1e+06 & 2.0e+06 & 0.0359, \ \ 0.0129
\end{tabular}
\end{center}
\centerline{{\bf Table 1.} (CMU classification errors depending on $\gamma$
  for $d=2$.)}
    
\bigskip

The number of cycles is the same as the number of
Cholesky factorizations that were
computed. The number of iterations is the number of steps in up- or down-
cycles where each iteration used at most ``order $n^2$ operations'',
somewhat less when the size
of the inactive set was small. (Since $n=500$ these examples used
less than $3n$ iterations.) The time is on a ThinkPad from 2016,
Intel(R) Core(TM) i7-6600U CPU @ 2.60GHz.
The KKT violation is the relative violation as detailed in
Section~\ref{sec.technical}.
The relative classification error is the relative number of points
which should have been classified in $S_1$ but were not, 
and likewise for $S_2$.

\subsubsection{Comparison with SMO}
In Section \ref{sec.param}
the overall classification errors for $\gamma = 0.03$ are slightly
lower than for the other values.
When choosing even smaller values of $\gamma = 0.03$ for this
example, the Hessian of $H$ is so poorly conditioned that a
reliable solution of \mr{qp} with the standard numerical
precision of about 16 decimal digits was not possible with any of
the methods.
For the value  $\gamma = 0.03$ the performance
of CMU was compared to SMO with greedy selection of the search
direction (GSMO) and with the much cheaper random selection (RSMO)
in the next table. To compensate for the cheaper iterations,
in GSMO the maximum number of iterations was set to 1000$n$
while it was set to 10000$n$ for RSMO (both, in Table 2 and Table 3).

\begin{center}
\begin{tabular}{ c c c c c c c }
  method & time & KKT violation &
  $q(x^{final})$ &  rel class. errors\\
  CMU  &  0.69  &   1.8e-11  &   -1.0e+12  & 0.0106, \ \ 0.0267\\
  GSMO & 17.9 & 2.0e-08 & -1.5e+10 & 0.3377, \ \   0.3537 \\
  RSMO & 22.0 & 1.8e-07 & -1.7e+09 & 0.3649, \ \   0.3051
\end{tabular}
\end{center}
\centerline{{\bf Table 2.} (Different algorithms for $\gamma=0.03$ and
   $d=2$.)}
    
\bigskip

The Hessian matrix for this problem is very ill-conditioned.
In spite of the fact that a large number of iterations
was allowed leading to a solution of reasonably high accuracy,
the classification error of GSMO or RSMO for this problem
was considerably higher than for CMU. (The minimum value of $q$
must be less or equal to the value
returned by CMU, indicating that also the values of $q$ generated
by GSMO or RSMO are far from optimality.)

The classification error is used, for example, in cross validation
approaches for the parameters of the kernel, and thus it is
important that the numerical accuracy is sufficiently high not
to deteriorate the classification error.
In the artificial example above, the accuracy generated by
either variant of SMO was not sufficient. 

There are many other versions of SMO with different choices of the pivot
element, but they all share the property that high accuracy solutions
require very many steps of SMO.

With a slightly different setting of $d\in\{3,5\}$ dimensions,
a much larger number of $n=10000$ training data points was compared on
a small computer cluster with 8 sockets,
64 CPUs, AMD Opteron(tm) Processor 6282 SE.
%3519102 jarre     20   0   15.0g   3.0g 204756 S 276.4   1.2
%30:24.42 MATLAB                                                                
The dimension of $H$ was 10000 by 10000,
and for such dimensions, the computation times
depend more closely on the numerical effort and to a lesser
extent on the overhead caused by the fact that the Matlab program used
is based on an (uncompiled) interpreter.
(The overhead grows about linearly  with the dimension
and the computational effort grows more than quadratically.)
The results are given in Table~3.

\begin{center}
\begin{tabular}{ c c c c c c c c c }
  method & $d$ & $\gamma$ & time & KKT violation & $\|x^{final}\|_\infty$ &
  $q(x^{final})$ &  rel class. errors\\
  CMU  & 3 & 0.03 & 24491 & 3.5e-10 & 5.0e+10 & -8.6e+12 & 0.0035, \ \  0.0170 \\
  GSMO & 3 & 0.03 & 18698 & 4.9e-09 & 4.9e+10 & -3.2e+11 & 0.2565, \ \  0.2619 \\
  RSMO & 3 & 0.03 & 11610 & 1.5e-06 & 5.7e+06 & -3.8e+09 & 0.1818, \ \  0.2921 \\ 
  CMU  & 5 & 0.03 & 24187 & 3.6e-10 & 4.7e+10 & -8.4e+12 & 0.0068, \ \  0.0223 \\
  GSMO & 5 & 0.03 & 19195 & 4.1e-07 & 2.2e+08 & -1.9e+10 & 0.4113, \ \  0.4458 \\
  RSMO & 5 & 0.03 & 11633 & 2.0e-06 & 6.7e+06 & -2.5e+09 & 0.1673, \ \  0.3509 \\
  CMU  & 5 &    3 & 27156 & 1.2e-14 & 2.4e+05 & -9.1e+05 & 0.0679, \ \  0.0789 \\
  GSMO & 5 &    3 & 5959  & 8.5e-11 & 2.4e+05 & -9.1e+05 & 0.0671, \ \  0.0774 \\
  RSMO & 5 &    3 & 11656 & 3.4e-04 & 1.8e+04 & -5.7e+05 & 0.0541, \ \  0.0833 
\end{tabular}
\end{center}
\centerline{{\bf Table 3.} (Different algorithms for $\gamma\in\{0.03,3\}$ 
    and $n=10000$.)}
    
\bigskip
Again, these numbers indicate that a small value of $\gamma$ should at least
be considered in a cross validation approach, and that for small values of
$\gamma$ a high numerical accuracy is required, in order to reduce the
classification error.

The number of outer iterations in CMU in these  examples
was between 6 and 8 with $11460$ to $14054$ inner iterations 
%11594 and 11460 and 14054
while GSMO and RSMO used 1000$n$ and 1000$n$ iterations respectively.
Only for $\gamma = 3$, GSMO terminated early after 350$n$ iterations
because the stopping criterion -- which was always tested after
integer multiples of $n$ iterations -- was satisfied. 

For CMU, the size of the first set of inactive indices
was limited to 2000 in order to
start with a moderately cheap Cholesky factorization, possibly at the
expense of more cycles. (This limitation does not apply to Table 2,
where $n=500$.)

Comparing GSMO and RSMO, it is
remarkable, that a more elaborate choice of pivots leads to
substantially higher numerical accuracy
(i.e., lower values of $q$)
even when considering the numerical
effort and allowing for 10 times more random steps.

\subsubsection{Dependence on $d$}
The next example illustrates the impact of the dimension $d$ on the
performance of CMU with $n=500$ training points.

\begin{center}
\begin{tabular}{ c c c c c c c c c }
  $d$ & cycles & iterations & time & KKT violation &
  $q(x^{final})$ & $\|x^{final}\|_\infty$ & rel class. errors\\
  2 &  3  & 433 &  0.69 & 1.8e-11 & -1.0e+12 & 5.8e+10 & 0.0106, \ \ 0.0267\\
  3 &  3  & 518 &  0.77 & 1.1e-11 & -1.2e+12 & 8.4e+10 & 0.0250, \ \ 0.0669\\
  5 &  7  & 686 &  1.05 & 2.4e-12 & -1.5e+11 & 5.8e+10 & 0.1728, \ \ 0.2001\\
 10 &  5  & 499 &  0.88 & 7.1e-15 & -1.7e+09 & 7.6e+08 & 0.4260, \ \ 0.3555\\
 50 &  3  & 260 &  0.60 & 9.8e-15 & -1.6e+07 & 6.8e+06 & 0.4849, \ \ 0.3786
\end{tabular}
\end{center}
\centerline{{\bf Table 4.} (Classification errors depending on $d$ for
  $\gamma=0.03$.)}
    
\bigskip

The condition numbers of $H\in\re^{n\times n}$ improve for larger values of $d$ but the
classification errors deteriorate.
(For $d>3$ more than $n=500$ points might be necessary to generate a classifier
with less than 10\% relative classification error.)
In any case, the purpose of
Table~4 was to illustrate the effect of higher dimensions of the data space 
on the running times and the numerical accuracy of CMU.

\subsection{\bf Checker board pattern}

In this example,
the sets $S_1$ and $S_2$ were defined along a $3\times 3$ checker board
pattern, and again 500 (uniformly randomly defined)
training points were used without errors
in the classifications of the training set.
Figure 1 shows the regions separated by the SVM with Gaussian kernel
with two different values of $\gamma$. The larger value of $\gamma$
results in ``more curvature'' of the boundary of the classification sets
(in green and red).

\vbox{
\vskip-4cm
\centerline{
  \includegraphics[width = 12cm]{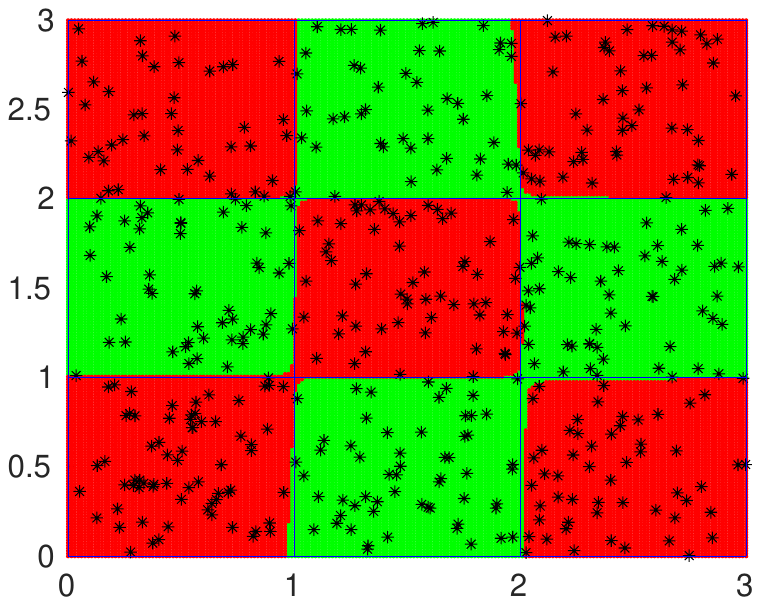} \hskip-3cm
  \includegraphics[width = 12cm]{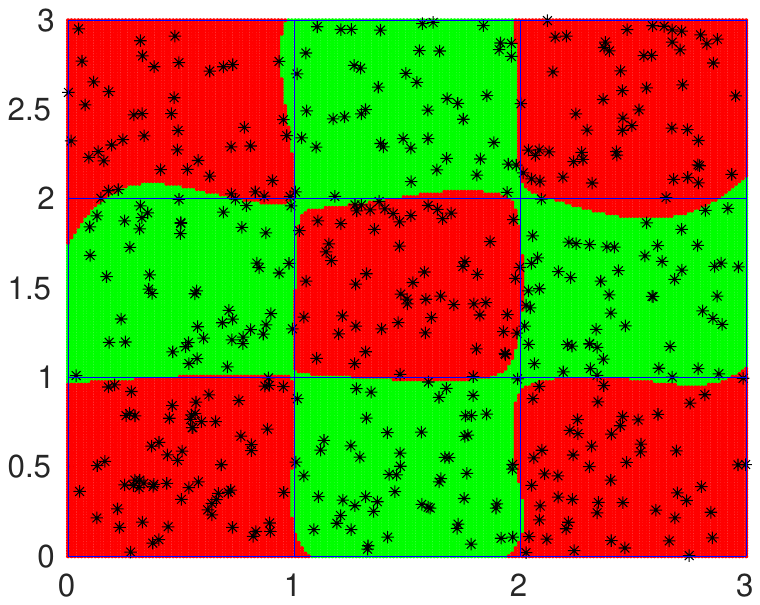}}
\vskip-5cm
\centerline{{\bf Figure 1}, CMU for Gaussian kernel with $C=\infty$ and
  $\gamma=0.03$ and $\gamma=3$ from left to right.}
\medskip
}

Figure 1 suggests that
the curvature of the boundary is lower for smaller values of $c$,
a fact that has been analyzed theoretically in \cite{jarresvm}
based on a modified self-concordance property of \cite{nene}.

\ignore{
A positive definite kernel induces a nonlinear mapping $\Phi$ of the data
points $x^{(i)}\in\re^d$ to a possibly higher-dimensional space $\Phi(\re^d)$,
and a linear separation of the image points $\Phi(x^i)$.
When $\Phi$ is an invertible affine mapping the linear separation
of the images  $\Phi(\re^d)$ translates to a linear separation
of the points  $x^{(i)}\in\re^d$, and when $\|\nabla^2\phi(x)\|$ is small, 
...
}

We point out that also in the plot on the right all training
points are classified correctly, and in fact, if no further information
is given, it might well be that
the pattern on the right describes the ``true pattern''.
For the ``exact checker board pattern'' that was actually used,
the above plot gives a visual indication that in some cases,
small values of the parameter $\gamma$ may be appropriate.

On the one side, as argued in \cite{jarresvm}, small values of $\gamma$
imply small values of a certain self-concordance parameter
similar to the one introduced in \cite{nene},
on the other side, small values of $\gamma$ result in a very
ill-conditioned matrix $H$. For $\gamma=0.03$ there were some
steps (less than 1 percent of the Newton steps),
where the Newton step resulted in a numerical increase
of the value of $q$. (We tested some of these instances: even though
the Newton step
$s$ satisfied $\nabla q(x)^Ts<-10^{-5}\|\nabla q(x) \|_2 \|s \|_2$
the computed values of $q$ did satisfy $q(x+\frac1{10}s)>q(x)$.) 
Therefore smaller values of $\gamma$ are not included in the
numerical examples of this section.

\bigskip

The same pattern with $\gamma=0.03$ is now used with GSMO
with 1000$n$ and 10000$n$ iterations. (Here, as well, RSMO
generated less accurate results.)

\vbox{
\vskip-4cm
\centerline{
  \includegraphics[width = 12cm]{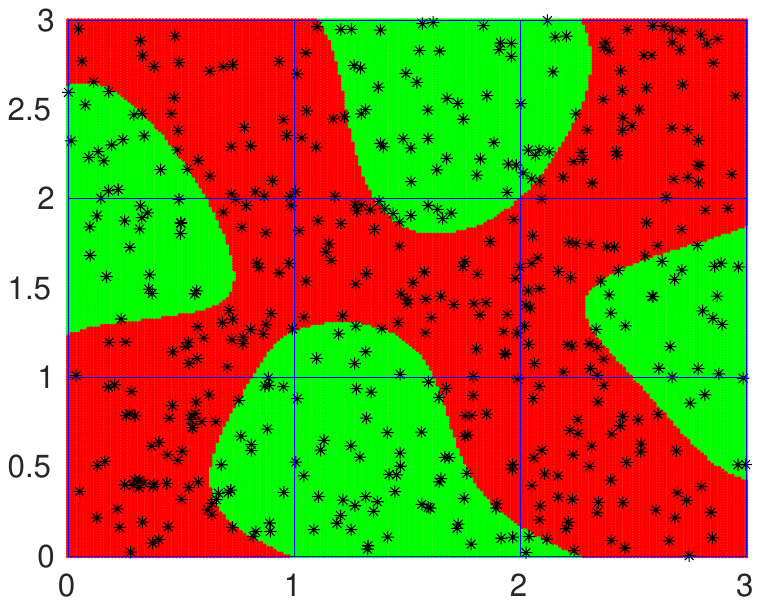} \hskip-3cm
  \includegraphics[width = 12cm]{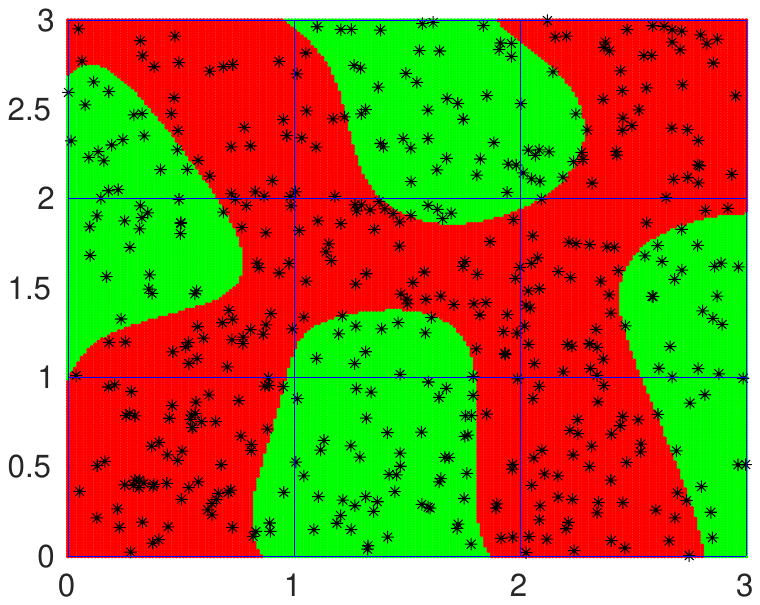}}
\vskip-5cm
\centerline{{\bf Figure 2}, Same example, illustrating that
  high numerical accuracy is essential.}
\centerline{GSMO for Gaussian kernel with 
  $\gamma=0.03$, 1000$n$ and 10000$n$ iterations from left to right.}
\medskip
}

Even when $10000n=10^8$ iterations are allowed, the accuracy of the solution
is rather low, and there are many falsely classified training data points.
Some of the results of above test runs are listed below. While the
tendency observed in Section \ref{sec.half} continues here as well
these examples are not intended to make a general claim about
the efficiency but to point out that there are instances for
which the computation of a high accuracy solution with
moderate computational effort may be essential.

\begin{center}
\begin{tabular}{ c c c c c c c c c }
  method &  iterations & time & KKT violation & 
  $q(x^{final})$ & $\|x^{final}\|_\infty$ \\
  CMU &  588    &  0.99 & 2.2e-11 & -1.4e+11 & 4.6e+10      \\
 GSMO & 500000  & 18.25 & 5.3e-08 & -6.2e+08 & 1.1e+08      \\
 GSMO & 5000000 & 197.7 & 9.5e-09 & -3.8e+09 & 5.7e+08
\end{tabular}
\end{center}
\centerline{{\bf Table 5.} (Checker board example, convergence of CMU and GSMO
  for $\gamma=0.03$ and $n=500$.)}
    
\bigskip

All files used to generate above data are available at

https://github.com/florianjarre/SVM-Test-Set

\section{Conclusion}
A key observation used in the CMU algorithm lies in the fact that repeated
stable and 
numerically cheap increases of the inactive set are possible while reducing
the objective function and without relying on a Hessian factorization. 
In the preliminary numerical experiments only a small number of cycles
were needed so that the overall numerical effort was small.
%while
%generating an exact optimal solution when there are no numerical rounding
%errors.
The implementation
in \cite{jarregithub}
includes a simple iterative refinement step that
helps reducing the numerical rounding errors.

\subsection*{Acknowledgment}
The author likes to thank Kevin Wischnewski for helpful comments
correcting and improving the presentation of this paper.

\end{document}